\documentclass[11pt]{article}
\usepackage{amsmath}
\usepackage{amsfonts}
\usepackage{amssymb}
\usepackage{enumerate}
\usepackage{amsthm}
\usepackage{graphicx}
\usepackage{verbatim}
\newtheorem{myproposition}{Proposition}[section]
\newtheorem{mytheorem}[myproposition]{Theorem}
\newtheorem{mylemma}[myproposition]{Lemma}

\newtheorem{mydefinition}[myproposition]{Definition}
\newtheorem{myalgorithm}[myproposition]{Algorithm}

\newtheorem{myexample}[myproposition]{Example}

\begin {document}

\title{Irregular Labellings of Circulant Graphs}
\author{Marcin Anholcer\thanks{E-mail:
m.anholcer@ue.poznan.pl}\\
{\small Pozna\'n University of Economics, Al. Niepodleg\l o\'sci 10, 60-967 Pozna\'n, Poland}\\
{\small Faculty of Informatics and Electronic Economy}\\
{\small Department of Operations Research}}
 \maketitle

\begin{abstract}
We investigate the \textit{irregularity strength} ($s(G)$) and \textit{total vertex irregularity strength} ($tvs(G)$) of circulant graphs $Ci_n(1,2,\dots,k)$ and prove that $tvs(Ci_n(1,2,\dots,k))=\left\lceil\frac{n+2k}{2k+1}\right\rceil$, while $s(Ci_n(1,2,\dots,k))=\left\lceil\frac{n+2k-1}{2k}\right\rceil$ except the case when $(n \bmod 4k = 2k+1 \wedge k\bmod 2=1) \vee n=2k+1$ and $s(Ci_n(1,2,\dots,k))=\left\lceil\frac{n+2k-1}{2k}\right\rceil+1$.
\end{abstract}


\section{Introduction}

Let us consider a simple undirected graph $G=(V(G),E(G))$ without loops, without isolated edges and with at most one isolated vertex. We assign a label $w(e)$ (called also weight), being natural positive number, to every edge $e \in E(G)$. For every vertex $v \in V(G)$ we define its \textit{weighted degree} as

\begin{eqnarray*}
wd(v)=\sum_{e\ni v}w(e).
\end{eqnarray*}

We call weighting $w$ \textit{irregular} if for each pair of vertices, their weighted degrees are distinct. In \cite{ref_ChaJacLehOelRuiSab1} the authors defined the graph parameter $s(G)$ called the irregularity strength of $G$ being the smallest integer $s$ such that there exists a weighting of $G$ with integers $\{1,2,\dots,s\}$ that is irregular. The value of $s(G)$ is known only for some special classes of graphs, e.g. complete graphs (\cite{ref_ChaJacLehOelRuiSab1}), graphs with the components being paths and cycles (\cite{ref_KinLeh},\cite{ref_AigTri}), or some families of trees (\cite{ref_Tog1}, \cite{ref_AmaTog}).

The lower bound on the $s(G)$ is given by the inequality

\begin{eqnarray}
s(G)\geq \max_{1\leq i\leq \Delta}\frac{n_i+i-1}{i} \label{ogr_dolne_ogolne}.
\end{eqnarray}

In the case of $d$-regular graphs it reduces to

\begin{equation}
s(G)\geq \frac{n+d-1}{d} \label{ogr_dolne_regularne}.
\end{equation}

The conjecture stated in \cite{ref_ChaJacLehOelRuiSab1} says that the value of $s(G)$ is for every graph equal to the above lower bound plus some constant not depending on $G$. The first upper bounds including the vertex degrees in the denominator were given in \cite{ref_FriGouKarPfe1} ($cn/\delta$ with relatively large values of $c$, depending on the relation between $n$, $\delta$ and $\Delta$), then improved (slight reduction of $c$) in \cite{ref_Prz1} and \cite{ref_Prz2}. The best upper bounds known so far can be found in \cite{ref_KalKarPfe1}. Namely, the authors have proved that

\begin{equation}
s(G)\leq \left\lceil\frac{6n}{\delta}\right\rceil \label{is_best}.
\end{equation}

The following variant of irregularity strength that allows also the vertices to be labeled was introduced in \cite{ref_BacJenMilRya1}. This time, the weighted degree is defined as 

\begin{eqnarray*}
wd(v)=\sum_{e\ni v}w(e)+w(v).
\end{eqnarray*}

The respective graph parameter, $tvs(G)$, is called \textit{total vertex irregularity strength}. The authors of \cite{ref_BacJenMilRya1} gave the following lower and upper bounds:

\begin{eqnarray}\label{Jendrol_bound1}
\left\lceil \frac{n+\delta(G)}{\Delta(G)+1} \right\rceil \leq tvs(G) \leq n+\Delta(G)-2\delta(G)+1.
\end{eqnarray}

In the case of $d$ - regular graphs it reduces to

\begin{eqnarray}\label{Jendrol_bound2}
\left\lceil \frac{n+d}{d+1} \right\rceil \leq tvs(G) \leq n-d+1.
\end{eqnarray}

The exact values of $tvs(G)$ are known only for few families of graphs, e.g. complete graphs, paths and cycles (\cite{ref_BacJenMilRya1}) or forests without vertices of degree $2$ (\cite{ref_AnhKarPfe}). The best upper bound on $tvs(G)$ is given in \cite{ref_AnhKalPrz1}:

\begin{equation}
tvs(G)\leq \left\lceil\frac{3n}{\delta}\right\rceil+1 \label{tvs_best}.
\end{equation}

Let us consider circulant graphs defined as follows (see e.g. \cite{ref_BarKheTog1}).

\begin{mydefinition}
Let $n$ and $s_1,s_2,\dots,s_k$ be integers, with $1\leq s_1<\dots<s_k\leq n/2$. The circulant graph $G=Ci_n(s_1,\dots,s_k)$ of order $n$ is a graph with vertex set $V(G)=\{0,1,\dots,n-1\}$ and edge set $E(G)=\{(x,x\pm s_i \mod n),x\in V(G), 1\leq i \leq k\}$.
\end{mydefinition}

Note that $Ci_n(s_1,\dots,s_k)$ is $2k$-regular. The main result given in \cite{ref_BarKheTog1} says that in the case $k=2$ and $s_1=1$,

\begin{equation}
s(Ci_n(1,s_2))=\left\lceil\frac{n+3}{4}\right\rceil
\end{equation}

\noindent{}if only $s_2\geq 2$ and $n\geq 4s_2+1$. Observe that in this case the value $s(G)$ is equal to the lower bound given by (\ref{ogr_dolne_regularne}).

In \cite{ref_AhmBac} in turn the authors gave the exact value of total vertex irregularity strength of the graphs $Ci_n(1,2)$:

\begin{equation}
tvs(Ci_n(1,2))=\left\lceil\frac{n+4}{5}\right\rceil.
\end{equation}

In this paper we consider more general case of circulant graphs, $Ci_n(1,2,\dots,k)$, i.e. the $k$-th powers of cycles $C_n^k$. We prove two following theorems.

\begin{mytheorem}\label{theorem_tvs}{\ }
\noindent{}If $k\geq 2$ and $n\geq 2k+1$, then
\begin{equation*}
tvs(C_n^k)=\left\lceil\frac{n+2k}{2k+1}\right\rceil.
\end{equation*}
\end{mytheorem}

\begin{mytheorem}\label{theorem_s}{\ }
\noindent{}If $k\geq 2$ and $n\geq 2k+1$, then
\begin{eqnarray*}
s(C_n^k))=\left\{
\begin{array}[]{ll}
\left\lceil\frac{n+2k-1}{2k}\right\rceil+1, & (n \bmod 4k = 2k+1 \wedge k\bmod 2=1) \vee n=2k+1,\\
\left\lceil\frac{n+2k-1}{2k}\right\rceil, &\text{otherwise.}
\end{array}
\right.
\end{eqnarray*}
\end{mytheorem}

\section{Proof of Theorem \ref{theorem_tvs}}

As the given value of ${ tvs}(G)$ equals to the lower bound given by (\ref{Jendrol_bound2}), it suffices to present the irregular weighting using the weights $1,\dots,s=\lceil\frac{n+2k}{2k+1}\rceil$.

Very briefly, the weighting of $G=C_n^k$ proceeds as follows.

We split $G$ into at most $s-1$ segments and label their edges in such a way, that the weighted degrees of vertices in every segment become distinct multiples of $2$. In the next step we multiply all the edge labels by about $s/2$ (depending on the parity of $s$) in order to obtain the weighting, where all the weighted degrees in any segment differ by at least $s-1$. By assigning distinct numbers from the set $\{1,2,\dots,s-1\}$ to the vertices in consecutive segments we obtain the final irregular weighting.

Now we are going to present the details of the proof. Let us start with two technical lemmas.

\begin{mylemma}\label{lemat_pojedyncze}{\ }

\noindent{}Let $S=S^{(k)}$ be graph on $2k+1$ vertices $\{v_0,v_1,\dots,v_{2k}\}$ (${k\geq 1}$) with  edge set consisting of the pairs $(v_i,v_{i+j})$, where ${i=0,1,\dots,2k-1}$ and~${j=1,2,\dots,\min\{k,2k-i\}}$. Assume there is a label $l(v_i)$ assigned to every vertex $v_i$, where
\begin{eqnarray*}
l(v_i)=
\begin{cases}
0 &\text{for } i\leq k,\\
2(i-k) &\text{for } i>k.
\end{cases}
\end{eqnarray*}

Then there exists labelling $w:E(S)\rightarrow \{0,1,2\}$ such that:

\begin{enumerate}[(i)]
\item
For every vertex $v_i$, $0\leq i \leq 2k$:
\begin{eqnarray*}
\sum_{e\ni v_i}w(e)+l(v_i)=2i.
\end{eqnarray*}
\item
Subgraph of $S$ consisting of the edges labeled $1$ contains all the vertices except $v_0$ and $v_{2k}$ and possesses an Euler cycle.
\end{enumerate}
\end{mylemma}

\noindent{}\textbf{Proof.} In order to prove the lemma we are going to present the algorithm that produces the desired weighting. We call \textit{closed} the vertices $v_i$ with all the edges $(v_i,v_j)$, $j>i$, weighted. The remaining vertices are called \textit{open}.

\begin{myalgorithm}\label{algorytm1}{\ }
\begin{enumerate}[(i)]
\item
Assign label $0$ to every edge $(v_0,v_i)$, where $1\leq i \leq k$. The vertex $v_0$ is now closed. Proceed to step $(ii)$.
\item
Let $i$ be the lowest integer such that all the vertices $v_j,0\leq j\leq i-1$, are closed and let
\begin{eqnarray*}
p=2i-\sum_{j=\max\{0,i-k\}}^{i-1}w((v_j,v_i))-l(v_i).
\end{eqnarray*}
If $1\leq p\leq k$, label forward edges according to the formula
\begin{eqnarray*}
w((v_i,v_j))=
\begin{cases}
0 &\text{for } i+1\leq j \leq i+k-p\\
1 &\text{for } i+k-p+1\leq j \leq i+k.
\end{cases}
\end{eqnarray*}
If $p\geq k+1$, use the formula
\begin{eqnarray*}
w((v_i,v_j))=
\begin{cases}
1 &\text{for } i+1\leq j \leq i+2k-p\\
2 &\text{for } j\geq i+2k-p+1.
\end{cases}
\end{eqnarray*}
Go to step $(iii)$.
\item
Close the vertex $v_i$. If $i=2k+1$, STOP. Otherwise go back to the step $(ii)$.
\end{enumerate}
\end{myalgorithm}

The number $p$ from step $(ii)$ always fulfils the condition $p\geq 1$ (see below).

Let us analyse the algorithm, considering its phases. For every vertex $v_i$, we call the edges $(v_i,v_j)$ \textit{forward} when $j>i$ and \textit{backward} when $j<i$. In every phase we process some group of vertices, by weighting its forward edges. Thus after each phase weighted degrees of vertices from some group reach their final values. 
 
For $0\leq i \leq \lfloor \frac{k}{2} \rfloor$, we assign the labels
\begin{eqnarray*}
w((v_i,v_j))=
\begin{cases}
0 &\text{for } i+1\leq j \leq k-i,\\
1 &\text{for } k-i+1\leq j\leq k+i.
\end{cases}
\end{eqnarray*}
After this phase the weighted degrees of vertices $v_i, i\leq \lfloor \frac{k}{2} \rfloor$, reach their final values $wd(v_i)=2i$. The weighted degrees of the remaining vertices obtain temporary values
\begin{eqnarray*}
wd(v_i)=
\begin{cases}
i-\lceil \frac{k}{2} \rceil &\text{for } \lfloor \frac{k}{2} \rfloor +1 \leq i \leq k,\\
i-\lceil \frac{k}{2} \rceil+1 &\text{for } k+1 \leq i \leq \lfloor \frac{3k}{2} \rfloor,\\
2(i-k) &\text{for } \lfloor \frac{3k}{2} \rfloor+1\leq i \leq 2k.
\end{cases}
\end{eqnarray*}
For $\lfloor \frac{k}{2} \rfloor+1\leq i \leq k$, the edges obtain the labels
\begin{eqnarray*}
w((v_i,v_j))=
\begin{cases}
1 &\text{for } i+1\leq j \leq k+i-1,\\
2 &\text{for } j=k+i.
\end{cases}
\end{eqnarray*}
After this phase the weighted degrees of the vertices $v_i, i \leq k$, reach their final values $wd(v_i)=2i$ and the remaining ones - temporary values
\begin{eqnarray*}
wd(v_i)=
\begin{cases}
i+1 &\text{for } k+1 \leq i \leq \lfloor \frac{3k}{2} \rfloor,\\
i+2 &\text{for } \lfloor \frac{3k}{2} \rfloor+1 \leq i \leq 2k.
\end{cases}
\end{eqnarray*}
In the next phase the edges incident to the vertices $v_i$, ${k+1\leq i \leq \lfloor \frac{3k}{2} \rfloor}$, obtain the labels
\begin{eqnarray*}
w((v_i,v_j))=
\begin{cases}
1 &\text{for } i+1\leq j \leq 3k-i,\\
2 &\text{for } 3k-i+1\leq j \leq 2k.
\end{cases}
\end{eqnarray*}
This way the weighted degrees of the vertices $v_i, i \leq \lfloor \frac{3k}{2} \rfloor$, reach their final values, and the remaining ones - temporary values
$$wd(v_i)=2i-k+1.$$
In the last phase all the remaining edges obtain label $2$. This way all the weighted degrees reach the final values
$$wd(v_i)=2i.$$
In order to prove the second part of the lemma, observe that after all the vertices having been processed, each of them except $v_0$ i $v_{2k}$ is incident to at least two edges labelled $1$. As all the weighted degrees are even, the degrees in the subgraph labeled $1$ are also even. Moreover this subgraph is connected as
$$w(v_i,v_k)=1,$$
for $1\leq i\leq 2k-1, i\neq k$.\qed

Second lemma guarantees the existence of analogous weighting for the segments of length $4k+2$.

\begin{mylemma}\label{lemat_podwojne}{\ }

\noindent{}Let $D=D^{(k)}$ be graph on $4k+2$ vertices $\{v_0,v_1,\dots,v_{4k+1}\}$ ($k\geq 2$) with edge set consisting of the pairs $(v_i,v_{i+j})$, where $i=0,1,\dots,4k$ and $j=1,2$, $\dots, {\min\{k,4k-i+1\}}$. Assume there is a label $l(v_i)$ assigned to every vertex $v_i$, where
\begin{eqnarray*}
l(v_i)=
\begin{cases}
0 &\text{for } i\leq 3k+1\\
2(i-3k-1) &\text{for } i>3k+1.
\end{cases}
\end{eqnarray*}
Then there exists labelling $w:E(D)\rightarrow \{0,1,2\}$ such that:
\begin{enumerate}[(i)]
\item
For every vertex $v_{i}$, $0\leq i \leq 4k+1$:
\begin{eqnarray*}
\sum_{e\ni v_{i}}w(e)+l(v_{i})=2\left\lfloor\frac{i}{2}\right\rfloor.
\end{eqnarray*}
\item
Subgraph of $D$ consisting of the edges labelled $1$ contains all the vertices except $v_0$, $v_1$, $v_{4k}$ and $v_{4k+1}$ and consists of at most two components, each of which possesses an Euler cycle.
\end{enumerate}
\end{mylemma}

\noindent{}\textbf{Proof.}

\noindent{}\textit{Case $1$: $k$ is even.}

\noindent{}Let us analyse the following algorithm, that produces the desired weighting. The \textit{open} and \textit{closed} vertices and the \textit{forward} and \textit{backward} edges are defined in the same way as in the case of algorithm \ref{algorytm1}\\

\begin{myalgorithm}\label{algorytm2}{\ }
\begin{enumerate}[(i)]
\item
Assign label $0$ to every edge $(v_0,v_i)$, where $1\leq i \leq k$. Vertex $v_0$ is closed now. Go to the step $(ii)$.
\item
Let $i$ be the lowest integer such that the vertices  $v_j$, ${0\leq j\leq i-1}$, are closed and let
\begin{eqnarray*}
p=2\lfloor i/2\rfloor -\sum_{i-k\leq j< i}w((v_j,v_i))-l(v_i).
\end{eqnarray*}
If $1\leq p\leq k$ label the forward edges according to the formula
\begin{eqnarray*}
w((v_i,v_j))=
\begin{cases}
0 &\text{for } i+1\leq j \leq i+k-p\\
1 &\text{for } i+k-p+1\leq j \leq i+k.
\end{cases}
\end{eqnarray*}
If $p\geq k+1$, use the formula
\begin{eqnarray*}
w((v_i,v_j))=
\begin{cases}
1 &\text{for } i+1\leq j \leq i+2k-p\\
2 &\text{for } i+2k-p+1\leq j \leq i+k.
\end{cases}
\end{eqnarray*}
Go to step $(iii)$.
\item
Close the vertex $v_i$. If $i=4k+2$, STOP. Otherwise go back to step $(ii)$.
\end{enumerate}
\end{myalgorithm}

The number $p$ from step $(ii)$ always fulfils the condition $p\geq 1$ (see below).

Similarly as in the case of Lemma \ref{lemat_pojedyncze}, we are going to analyse the algorithm phase by phase, each time considering chosen set of vertices and its forward edges.

For $0 \leq i \leq \frac{k-2}{2}$, the forward edges obtain labels
\begin{eqnarray*}
w((v_{2i},v_j))=
\begin{cases}
0 &\text{for } 2i+1\leq j \leq k,\\
1 &\text{for } k+1\leq j \leq 2i+k,
\end{cases}
\end{eqnarray*}
and
\begin{eqnarray*}
w((v_{2i+1},v_j))=
\begin{cases}
0 &\text{for } 2i+2\leq j \leq k+1,\\
1 &\text{for } k+2\leq j \leq 2i+k+1.
\end{cases}
\end{eqnarray*}
Moreover
$$w((v_k,v_j))=1,{\ \ \ } k+1\leq j \leq 2k,$$
and
\begin{eqnarray*}
w((v_{k+1},v_j))=
\begin{cases}
0 &\text{for } k+2\leq j \leq \frac{3k}{2}+1,\\
1 &\text{for } \frac{3k}{2}+2\leq j \leq 2k+1.
\end{cases}
\end{eqnarray*}
After this phase weighted degrees of the vertices $v_i, i \leq k+1$ reach their final values 
$$wd(v_i)=2\left\lfloor\frac{i}{2}\right\rfloor,$$
while the remaining ones - temporary values
\begin{eqnarray*}
wd(v_i)=
\begin{cases}
2k-i+1 &\text{for } k+2\leq i \leq \frac{3k}{2}+1,\\
2k-i+2 &\text{for } \frac{3k}{2}+2 \leq i \leq 2k+1,\\
0 &\text{for } 2k+2 \leq i \leq 3k+1,\\
2(i-3k-1) &\text{for } 3k+2 \leq i \leq 4k+1.
\end{cases}
\end{eqnarray*}

For $k+2\leq i \leq \frac{3k}{2}$, the forward edges obtain labels
\begin{eqnarray*}
w((v_i,v_j))=
\begin{cases}
0 &\text{for } i+1\leq j \leq i+k-3\lceil \frac{i-k-1}{2}\rceil-\lfloor\frac{i-k-1}{2}\rfloor,\\
1 &\text{for } i+k-3\lceil \frac{i-k-1}{2}\rceil-\lfloor\frac{i-k-1}{2}\rfloor+1\leq j \leq i+k.
\end{cases}
\end{eqnarray*}
After this phase weighted degrees of all the vertices $v_i$, $0\leq i \leq \frac{3k}{2}$ reach the final values, and the remaining weighted degrees - temporary
\begin{eqnarray*}
wd(v_i)=
\begin{cases}
\frac{k}{2} &\text{for } i=\frac{3k}{2}+1,\\
\frac{k}{2}+\lfloor\frac{i}{2}\rfloor-\lceil\frac{i}{2}\rceil+1 &\text{for } \frac{3k}{2}+2 \leq i \leq 2k+1,\\
\frac{5k}{2}-i+1 &\text{for } 2k+2 \leq i \leq \frac{5k}{2},\\
0 &\text{for } \frac{5k}{2}+1 \leq i \leq 3k+1,\\
2(i-3k-1) &\text{for } 3k+2 \leq i \leq 4k+1.
\end{cases}
\end{eqnarray*}

In the next phase we put $1$ on the edges $(v_{\frac{3k}{2}+1},v_j)$ for ${\frac{3k}{2}+2\leq j\leq \frac{5k}{2}}$. If $\frac{k}{2}$ is odd, the edge $(v_{\frac{3k}{2}+1},v_{\frac{5k}{2}+1})$ obtains label $1$, otherwise $2$. For $\frac{3k}{2}+2 \leq i \leq 2k+1$ we assign the labels
\begin{eqnarray*}
w((v_i,v_j))=1, {\ \ \ }\text{for } i+1\leq j \leq i+k.
\end{eqnarray*}

The weighted degrees of the vertices $v_{i}$, $0\leq i \leq 2k+1$, reach now their final values, while the remaining ones are equal to
\begin{eqnarray*}
wd(v_i)=
\begin{cases}
3k+2-i &\text{for } 2k+2 \leq i \leq 3k+1,i\neq \frac{5k}{2}+1,\\
2\left\lceil\frac{k}{4}\right\rceil+1 &\text{for } i=\frac{5k}{2}+1,\\
2(i-3k-1) &\text{for } 3k+2 \leq i \leq 4k+1.
\end{cases}
\end{eqnarray*}

For $2k+2\leq i \leq 3k-1$, $i \neq \frac{5k}{2}+1$, the edges obtain labels
\begin{eqnarray*}
w((v_i,v_j))=
\begin{cases}
1 &\text{for } i+1\leq j \leq i+k-2\left\lfloor \frac{i-2k}{2} \right\rfloor,\\
2 &\text{for } i+k-2\lfloor \frac{i-2k}{2} \rfloor+1\leq j\leq i+k,
\end{cases}
\end{eqnarray*}
and for $i=\frac{5k}{2}+1$ - labels
\begin{eqnarray*}
w((v_i,v_j))=
\begin{cases}
1 &\text{for } i+1\leq j \leq 3k+1,\\
2 &\text{for } 3k+2\leq j\leq i+k.
\end{cases}
\end{eqnarray*}

After this assignment, the weighted degrees of the vertices $v_i$, $0\leq i \leq 3k-1$ reach their final values. The remaining ones are equal to
\begin{eqnarray*}
wd(v_i)=
\begin{cases}
k &\text{for } i=3k,\\
2\left\lfloor\frac{3k}{4}\right\rfloor-2 &\text{for } i=3k+1,\\
2k-2 &\text{for } 3k+2 \leq i \leq 4k,\\
2k &\text{for } i=4k+1.
\end{cases}
\end{eqnarray*}

For $i=3k$, the weights obtain values
$$w((v_{3k},v_j))=2,{\ \ \ }\text{for } 3k+1\leq j\leq 4k.$$
For $i=3k+1$ in turn we label forward edges with
\begin{eqnarray*}
w((v_i,v_j))=
\begin{cases}
1 &\text{for } i+1\leq j \leq 2\left\lfloor\frac{7k}{4}\right\rfloor+1,\\
2 &\text{for } 2\left\lfloor\frac{7k}{4}\right\rfloor+2\leq j\leq i+k.
\end{cases}
\end{eqnarray*}

After this phase the weighted degrees of vertices $v_i$, where $0\leq i \leq 3k+1$, reach their final values. The remaining ones are equal to
\begin{eqnarray*}
wd(v_i)=
\begin{cases}
2k+1 &\text{for } 3k+2 \leq i \leq \frac{7k}{2},\\
2\left\lceil\frac{k}{4}\right\rceil+\frac{3k}{2}+1 &\text{for } i=\frac{7k}{2}+1,\\
2k+2 &\text{for } \frac{7k}{2}+2 \leq i \leq 4k+1.
\end{cases}
\end{eqnarray*}

For $3k+2 \leq i \leq 2\lfloor\frac{7k}{4} \rfloor+1$, the edges are labelled with
\begin{eqnarray*}
w((v_i,v_j))=
\begin{cases}
1 &\text{for } i+1\leq j\leq 4k-2\lfloor\frac{i-3k}{2}\rfloor+1,\\
2 &\text{for } 4k-2\lfloor\frac{i-3k}{2}\rfloor+2\leq j\leq 4k+1.
\end{cases}
\end{eqnarray*}

After this assignment the weighted degrees of vertices $v_i$, where $0 \leq i \leq 2\lfloor\frac{7k}{4} \rfloor+1$, obtain its final values. The remaining ones are equal to
$$wd(v_i)=2\left\lfloor\frac{i}{2}\right\rfloor-k-2(\frac{k}{2}\bmod 2)+2.$$

In the last phase we assign label $2$ to all the remaining edges. This produces the final irregular weighting. 

\vspace{10pt}
\noindent{}\textit{Case $2$: $k$ is odd.}

\noindent{}Let us consider subgraph of $D=D^{(k)}$ isomorphic to $D^\star=D^{(k^\star)}$, $k^\star=k-1$, with the vertex set $V(D^\star)=(v^\star_0, v^\star_1, \dots, v^\star_{4k^\star+1})$, where $v^\star_i=v_{i-4}$ for ${i=0,\dots,4(k-1)+1}$. In order to label the edges of $D^\star$ we use the algorithm \ref{algorytm2}. Let us denote the initial vertex labels by $l^\star(v_i)$. Obviously $l^\star(v_i)=l(v_i)-2$ for $3k+2\leq i \leq 4k+1$. Let us consider partial weighting of $D$ obtained this way. Let $wd^\star(v_i)$ be temporary weighted degree of $v_i$, while $wd(v_i)$ its final value. Observe that
\begin{eqnarray*}
wd^\star(v_i)=\left\{
\begin{array}[]{ll}
wd(v_i)&\text{for } 0 \leq i \leq 1\\
wd(v_i)-2&\text{for } 2 \leq i \leq 3\\
wd(v_i)-4&\text{for } 4 \leq i \leq 4k+1
\end{array}
\right.
\end{eqnarray*}
Let us change the vertex labels $l^\star(v_i)$ to $l(v_i)$. This increases by $2$ the values of $wd^\star(v_i)$ for $3k+2 \leq i \leq 4k+1$.

The edges $(v_i, v_{i+k})$ for $4 \leq i \leq 3k+1$ have not been weighted yet. We assign label $2$ to them. As the weighted degrees increase by $2$ or $4$, depending on the number of incident edges, we obtain:
\begin{eqnarray*}
wd^\star(v_i)=\left\{
\begin{array}[]{ll}
wd(v_i)&\text{for } 0 \leq i \leq 1\\
wd(v_i)-2&\text{for } 2 \leq i \leq k+3\\
wd(v_i)&\text{for } k+4 \leq i \leq 4k+1
\end{array}
\right.
\end{eqnarray*}
If $k\geq 5$, then the edges $(v_i, v_{i+1})$ for $6\leq i \leq k+2$ have label $0$ (they correspond with the edges $(v^\star_i, v^\star_{i+1})$ of $D^\star$ for $2\leq i \leq k^\star-1$). The same is the label of the edge $(v_6,v_{k+3})$ (corresponding with $(v^\star_2, v^\star_{k^\star})$ in $D^\star$). We change them to $1$. Similarly, the edges $(v_2,v_3)$, $(v_3,v_4)$, $(v_4,v_5)$ and $(v_2,v_5)$ are labelled $0$, which we change to $1$. We put $0$ on all not labelled edges of $D$ (not belonging to $D^\star$).

If $k=3$, then we put $1$ on the edges $(v_2,v_3)$, $(v_3,v_6)$, $(v_5,v_6)$, $(v_4,v_5)$ and $(v_2,v_4)$ (all of which are either labelled $0$ or not labelled yet). This way we obtain the desired weighting of $D$.

To prove part $(ii)$ of the Lemma observe that after processing all the vertices, each of them except $v_0$, $v_1$, $v_{4k}$ and $v_{4k+1}$ is incident to at least two edges labelled $1$. As weighted degrees are even, the degrees in the subgraph labelled $1$ are also even. Moreover, the subgraph is connected when $k$ is even or $k=3$, otherwise it consists of two components. When $k$ is even, label $1$ is assigned in particular to the edges
\begin{eqnarray*}
\begin{array}[]{ll}
(v_i,v_{k+2})& {\ \ \ }\text{for } 2 \leq i \leq k,\\
(v_k,v_i)& {\ \ \ }\text{for } k+1 \leq i \leq 2k,\\
(v_{2k},v_i)& {\ \ \ }\text{for } 2k+1 \leq i \leq 3k,\\
(v_{2k+1},v_{3k+1}),& \\
(v_{3k+1},v_{3k+2}),& \\
(v_{3k+2},v_{i})& {\ \ \ }\text{for } 3k+3\leq i \leq 4k-1.
\end{array}
\end{eqnarray*}

In the case when $k\geq 5$ is odd, respective edges form an Euler graph with vertex set $v_i$, $6\leq i \leq 4k-1$ and the second component is subgraph on the vertices $v_2$, $v_3$, $v_4$ and $v_5$. If $k=3$, then the Euler graph is formed by the edges $(v_2,v_3)$, $(v_2,v_4)$, $(v_3,v_6)$, $(v_4,v_5)$, $(v_5,v_6)$, $(v_6,v_7)$, $(v_6,v_8)$, $(v_7,v_9)$, $(v_8,v_9)$, $(v_9,v_{10})$, $(v_9,v_{11})$ and $(v_{10},v_{11})$. \qed

Let us go back to the proof of Theorem \ref{theorem_tvs}. We are going to consider few cases, depending on the relation between $n$ and $k$.

\vspace{10pt}
\noindent{}\textit{Case $1$: $n=2k+1$.}

In this case $C_n^k$ is isomorphic to $K_n$ and the proof of the equality
\begin{eqnarray*}
tvs(C_n^k)=tvs(K_n)=2
\end{eqnarray*}
may be found e.g. in the paper of M. Ba\v ca et al. \cite{ref_BacJenMilRya1}.

\vspace{10pt}
\noindent{}\textit{Case $2$: $n=2k+2$.}

Let us assign to the vertices of $C_n^k$ indices $-k$, $-k+1$, $\dots$, $-1$, $0$, $1$, $2$, $\dots$, $k$, $k+1$. We assign the labels to the edges in the following way:
\begin{eqnarray*}
w(v_i,v_j)=\left\{
\begin{array}[]{ll}
1& \text{for } |i|+|j|\leq k\vee (|i|+|j|=k+1 \wedge \max{\{i,j\}}\leq 0), \\
2& \text{otherwise.}
\end{array}
\right.
\end{eqnarray*}
After such assignment, the weighted degrees are equal to:
\begin{eqnarray*}
wd(v_i)=\left\{
\begin{array}[]{ll}
2k-2i-2&\text{for } -k \leq i \leq -\lceil\frac{k}{2}\rceil-1 \\
2k-2i-1&\text{for } -\lceil\frac{k}{2}\rceil \leq i \leq -1 \\
2k+2i&\text{for } 0 \leq i \leq \lfloor\frac{k}{2}\rfloor \\
2k+2i-1&\text{for } \lfloor\frac{k}{2}\rfloor+1 \leq i\leq k \\
2k+2i-2&\text{for } i=k+1
\end{array}
\right.
\end{eqnarray*}

It means that the weighted degrees of vertices $v_i$, $-\lceil\frac{k}{2}\rceil\leq i \leq \lfloor\frac{k}{2}\rfloor$, are distinct integers from the set $\{2k,\dots,3k\}$, while the remaining ones - distinct integers from the set $\{3k,\dots,4k\}$.In order to obtain the final irregular weighting we assign $1$ to all the vertices from the first group, and $2$ to the remaining ones. 

\noindent{}\textit{Case $3$: $n>2k+2$.}

We can express the number of vertices of $G=C_n^k$ as ${n=t(4k+2)+r}$, where $t\geq 0$ and $1\leq r\leq 4k+2$ are some integers.

We split $G$ into some number of segments, each of which except two are isomorphic to $S^{(k)}$, one to $S^{(g)}$ and one to $S^{(h)}$ or $D^{(h)}$, where $g,h\leq k$, $|g-h|\leq 1$. For some values of $r$ it is necessary to include one additional vertex. Using Lemmas \ref{lemat_pojedyncze} and \ref{lemat_podwojne} we construct the weightings for the segments, and then we expand it to whole graph $G$.

There are some differences in the construction depending on the relation between $r$ and $k$, so we consider three cases. As some concepts and elements of the construction repeat, for convenience we consider first the case $2k+3\leq r\leq 4k+2$, then $2\leq r\leq 2k+1$ and finally $r\in\{1,2k+2\}$.

\vspace{10pt}
\noindent{}\textit{Case $3.1$: $2k+3\leq r\leq 4k+2$.}

Let $g=\left\lfloor\frac{r-2}{4}\right\rfloor$. Consider $2t$ graphs $S_1, S_2,\dots, S_{2t}$ isomorphic to $S^{(k)}$ and two additional graphs $S_{2t+1}\cong S^{(h)}$ and $S_{2t+2}\cong S^{(g)}$ where $h=g$ if $r$ is odd and $h=g+1$ otherwise. 

Let us denote the vertices of $G$ by $v_0(G),\dots,v_{n-1}(G)$. We split $G$ into segments in such a way that the vertices of the graphs $S_j$, $1\leq j\leq 2t$ (denoted by $v_0(S_j),\dots,v_{2k+1}(S_j)$), $S_{2t+1}$ (denoted by $v_0(S_{2t+1}),\dots,v_{2h}(S_{2t+1})$) and $S_{2t+2}$ (denoted by $v_0(S_{2t+2}),\dots,v_{2g}(S_{2t+2})$) are identified with consecutive subsets of vertices of $G$, where for increasing indices of vertices of $G$ the indices of vertices in $S_j$ increase for $j=2i-1$ and decrease for $j=2i$, $1\leq i\leq t+1$. In other words, if the vertices $v_i(G)$ and $v_{i+1}(G)$ (or $v_{n-1}(G)$ and $v_{0}(G)$) belong to two consecutive segments $S_j$ and $S_{j+1}$ (respectively $S_{2t+2}$ and $S_1$), then either they both have highest indices in the segments (so $v_i(G)=v_{2k}(S_j)$ and $v_{i+1}(G)=v_{2k}(S_{j+1})$ or possibly $v_i(G)=v_{2h}(S_j)$ and $v_{i+1}(G)=v_{2g}(S_{j+1})$) or $v_i(G)=v_{0}(S_j)$ and $v_{i+1}(G)=v_{0}(S_{j+1})$. If $r \bmod 4 \in \{1,3\}$, we include additional vertex $v^\star=v_{n-2h-2}(G)$ between $S_{2t+1}$ and $S_{2t+2}$.

We perform weighting of $G$ in two steps. We begin with the weighting of edges $w^\star:E(G)\rightarrow \{0,1,2\}$, then modify weights of edges and label the vertices in order to obtain irregular weighting $w:E(G)\cup V(G)\rightarrow \{1,\dots,s\}$.

We begin the construction of $w^\star$ with weighting the segments using Algorithm \ref{algorytm1}. Then we assign labels to the remaining edges of $G$ (i.e. edges $e\in E(G)\setminus \bigcup_{j=1}^{2t+2}E(S_j)$). 

If two consecutive segments are $S_{2i-1}$ and $S_{2i}$, $1\leq i\leq t$ (so two consecutive vertices belonging to them have weighted degree $2k$), we put $2$ on all the edges joining them (not belonging to the segments, but belonging to $G$). If the segments are $S_{2i}$ and $S_{2i+1}$, $1\leq i\leq t$ or $S_{2t+2}$ and $S_1$, then we assign $0$ to the edges. Observe that the weighted degrees do not change as the sums of labels of the edges joining different segments of $G$ are equal to the vertex labels $l(v)$ from Lemma \ref{lemat_pojedyncze}.

If $r \bmod 2 =0$, we put $2$ on the edges
\begin{eqnarray*}
\begin{array}[]{ll}
(v_{h+i}(S_{2t+1}),v_{2g-j+1}(S_{2t+2})) &\text{for }  1\leq i\leq h, 1\leq j\leq i,
\end{array}
\end{eqnarray*}
otherwise on the edges
\begin{eqnarray*}
\begin{array}[]{ll}
(v_{h+i}(S_{2t+1}),v_{2g-j+2}(S_{2t+2})) &\text{for }  2\leq i\leq h, 1\leq j\leq i-1,\\
(v^\star,v_{h+i}(S_{2t+1})) &\text{for }  1\leq i\leq g+1,\\
(v^\star,v_{2g-i+1}(S_{2t+2})) &\text{for }  1\leq i\leq g+1.
\end{array}
\end{eqnarray*}

Observe that last operation may increase by $2$ weighted degrees of some vertices $v\in V(S_{2t+2})$, however the conditions
$$wd^\star(v)\leq 2k {\ \ \ }\text{for }  v\in V(S_{2t+1})\cup V(S_{2t+2})\cup \{v^\star\}$$
and
$$|\{v\in V(S_{2t+1})\cup V(S_{2t+2})\cup \{v^\star\}: wd^\star(v)=2i\}|\leq 2{\ \ \ }\text{for } i=0,\dots,k,$$
remain satisfied as $wd(v^\star)\neq wd(v_{2g}(S_{2t+2}))$.

We finish the construction of $w^\star$ by putting $0$ on all the edges that have not been labelled so far $G$.

Finally we obtain the weighting $w^\star:E(G)\rightarrow \{0,1,2\}$ such that all the weighted degrees are even, the subgraph formed by the edges labelled $1$ has at most $2t+2$ components being Euler graphs and
$$|\{v\in V(G): wd^\star(v)=2i\}|\leq 2t+2{\ \ \ }\text{for } i=0,\dots,k.$$

Let $s=\lceil \frac{n+2k}{2k+1} \rceil$ (observe that $s$ is odd). We define weighting $w:E(G)\cup V(G)\rightarrow \{1,2,\dots, s\}$ in the following way:

\begin{eqnarray}\label{nowe_we}
w(e)=\frac{s-1}{2}w^\star(e)+1{\ \ \ }\text{for }  e\in E(G),
\end{eqnarray}
\begin{eqnarray}\label{nowe_wv}
w(v)=\left\{
\begin{array}[]{ll}
j&\text{for }  v\in V(S_j),1\leq j\leq 2t+1,\\
2t+2&\text{for }  v\in V(S_{2t+2})\cup \{v^\star\}.
\end{array}
\right. 
\end{eqnarray}

As we can see, in the modified weighting $w:E(G)\rightarrow \{0,\frac{s+1}{2},s\}$ the sums of labels of edges incident to distinct vertices of every segment $S_j$ differ by at least $s-1$. By putting $j$ on the vertices of $S_j$, we distinguish weighted degrees in whole graph $G$ as $2t+2\leq s-1$.

\vspace{10pt}
\noindent{}\textit{Case $3.2$: $2\leq r\leq 2k+1$.}

Let $g_1=\left\lfloor\frac{4k+r+1}{6}\right\rfloor$. If $g_1\geq 2$ then we proceed as in the previous case, this time using $2t-2$ graphs $S_1, S_2,\dots S_{2t-2}$ isomorphic to $S^{(k)}$ and two additional graphs $S_{2t}\cong D^{(g_1)}$ and $S_{2t-1}\cong S^{(g_2)}$, where
\begin{eqnarray*}
g_2=
\begin{cases}
g_1 &\text{for }  (4k+r \bmod 6) \in \{1,2\},\\
g_1+1 &\text{for }  (4k+r \bmod 6) \in \{3,4\},\\
g_1-1 &\text{for }  (4k+r \bmod 6) \in \{0,5\}.
\end{cases}
\end{eqnarray*}

We label all the graphs $S_j,1\leq j\leq 2t-1$ using the Algorithm \ref{algorytm1}, and the graph $S_{2t}$ using the Algorithm \ref{algorytm2}, possibly with additional modification when $g_1$ is odd (see proof of Lemma \ref{lemat_podwojne}).

\begin{figure}[h]
    \begin{center}
    	\scalebox{0.5}{\includegraphics{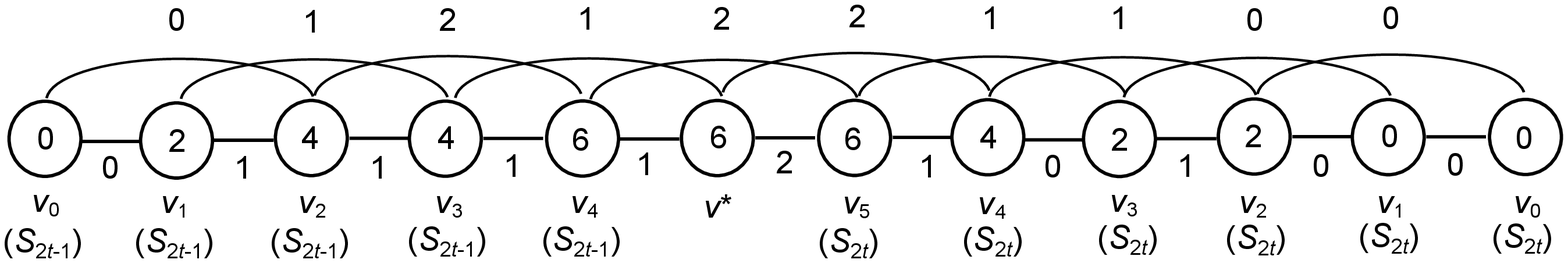}}
    \end{center}
    \caption{The labelling of $S_{2t-1}$ and $S_{2t}$ and the edges incident with $v^\star$ when $g_1=1$}
    \label{circulant_rys_specjalne_gD1}
\end{figure}

Let $g=\min\{g_1,g_2\}$ and $h=\max\{g_1,g_2\}$. If $g_1=1$ (then $g_2=2$), in the case of graphs $S_{2t-1}$ and $S_{2t}$ we use the weighting as in the Figure \ref{circulant_rys_specjalne_gD1} instead of the mentioned algorithms (we have to add vertex $v^\star$, and all the edges not included in the figure, either belonging to $E(S_{2t-1})\cup E(S_{2t})$ or joining these graphs one with other or with $v^\star$, are labelled $0$).

We identify the vertices of the segments $S_j,1\leq j \leq 2t$, with consecutive subsets of vertices of $G$ and label all remaining edges (except the ones joining $S_{2t-1}$ with $S_{2t}$) as in the previous case. This time we include additional vertex $v^\star$ between $S_{2t-1}$ and $S_{2t}$ if $r$ is even and $g_1>1$. In such case we put $2$ on the edges
\begin{eqnarray*}
\begin{array}[]{ll}
(v_{2g_2-h+i}(S_{2t-1}),v_{4g_1-j+3}(S_{2t})) &\text{for }  2\leq i\leq h, 1\leq j\leq i-1,\\
(v^\star,v_{2g_2-h+i}(S_{2t-1})) &\text{for }  1\leq i\leq g+1,\\
(v^\star,v_{4g_1-i+2}(S_{2t})) &\text{for }  1\leq i\leq g+1.
\end{array}
\end{eqnarray*}

If  $r$ is odd, we assign $2$ to the edges
\begin{eqnarray*}
\begin{array}[]{ll}
(v_{2g_2-h+i}(S_{2t-1}),v_{4g_1-j+2}(S_{2t})) &\text{for }  1\leq i\leq h, 1\leq j\leq i.
\end{array}
\end{eqnarray*}

If $g_1=1$ then we append $v^\star$ as in the Figure \ref{circulant_rys_specjalne_gD1}.

We put $0$ on all the remaining edges. Finally we obtain weighting $w^\star:E(G)\rightarrow \{0,1,2\}$ such that all weighted degrees are even, the subgraph labelled $1$ has at most $2t+1$ components being Euler graphs, and
$$|\{v\in V(G): wd^\star(v)=2i\}|\leq 2t+1{\ \ \ }\text{for } i=0,\dots,k.$$

This time $s$ is even, so we have to change slightly the construction of $w$. Firstly, we use the formula (\ref{nowe_we}) only for the edges $e$ with $w^\star(e)\in\{0,2\}$. In the case $w^\star(e)=1$ we proceed as follows. For every component labelled $1$ (being Euler graph) we start in any vertex and moving through an Euler walk label the edges alternately with:

\begin{eqnarray}\label{nowe_we2}
\begin{array}[]{ll}
w^{(1)}(e)=\frac{s-2}{2}w^\star(e)+1,\\
w^{(2)}(e)=\frac{s}{2}w^\star(e)+1.
\end{array}
\end{eqnarray}

If the walk has even length, new weighted degrees in every segment will differ by at least $s-1$, as in the previous case. If it is odd, the only exceptions are the starting vertices of the walks and their neighbours, as the starting ones may have the weighted degree lower by $1$ than the desired one. Let $V_{0}$ be the set of all such vertices. In order to finish the construction we assign to the vertices labels
\begin{eqnarray}\label{nowe_wv2}
w(v)=\left\{
\begin{array}[]{ll}
j&\text{for }  v\in V(S_j),1\leq j\leq 2t-2,\\
w_0\in \{2t-1,2t,2t+1\}&\text{for }  v\in V(S_{2t+2})\cup \{v^\star\}.
\end{array}
\right. 
\end{eqnarray}
The numbers $w_0\in \{2t-1,2t,2t+1\}$ have to be assigned in such a way, that every two vertices $v^\prime, v^{\prime\prime} \in V(S_{2t-1})\cup V(S_{2t})\cup\{v^\star\}$, for which
$$wd^\star(v^{\prime})=wd^\star(v^{\prime\prime}),$$
obtain two distinct labels. Such an assignment is possible as 
$$|\{v\in V(S_{2t-1})\cup V(S_{2t})\cup\{v^\star\}: wd^\star(v)=2i\}|\leq 3{\ \ \ }\text{for } i=0,\dots,k.$$

In order to finish the construction we increase by $1$ the label of every vertex $v\in V_0$. The obtained weighting is irregular as ${2t+1\leq s-1}$.

\vspace{10pt}
\noindent{}\textit{Case $3.3$: $r=1$ or $r=2k+2$.}

In such a situation we find irregular weighting for $C_{n-1}^k$ as in one of the previous cases. Note that the maximum weighted degree equals to $(2k+1)s-1$ (for some vertex labelled $s-1$ with all incident edges labelled $s$). Let $v_i$ and $v_{i+1}$ be two adjacent vertices with all incident edges labelled $s$. We remove the edges
\begin{eqnarray*}
\begin{array}[]{ll}
(v_j, v_{j+k}) &\text{for }  i-k+1\leq j \leq i,\\
(v_{j-k}, v_j) &\text{for }  i+1 \leq j \leq i+k,
\end{array}
\end{eqnarray*}
and append the vertex $v_{n-1}$ and the edges
\begin{eqnarray*}
\begin{array}[]{ll}
(v_j,v_{n-1}) &\text{for }  i-k+1\leq j \leq i,\\
(v_{n-1},v_j) &\text{for }  i+1 \leq j \leq i+k,
\end{array}
\end{eqnarray*}
all labelled $s$. Finally we put $w(v_{n-1})=s$.\qed

\begin{myexample}[Irregular total weighting of $C_{22}^2$]{\ }
\\ \noindent{}We have $n=22$, $k=2$, $s=6$, $t=4$, $r=2$, $g_1=1$, $g_2=2$, $g=1$ and $h=2$. Weighting $w^\star$ is presented on Figure \ref{circulant_przyklad_1}. The edges joining different segments are red.

\begin{figure}[h]
    \begin{center}
    	\scalebox{0.5}{\includegraphics{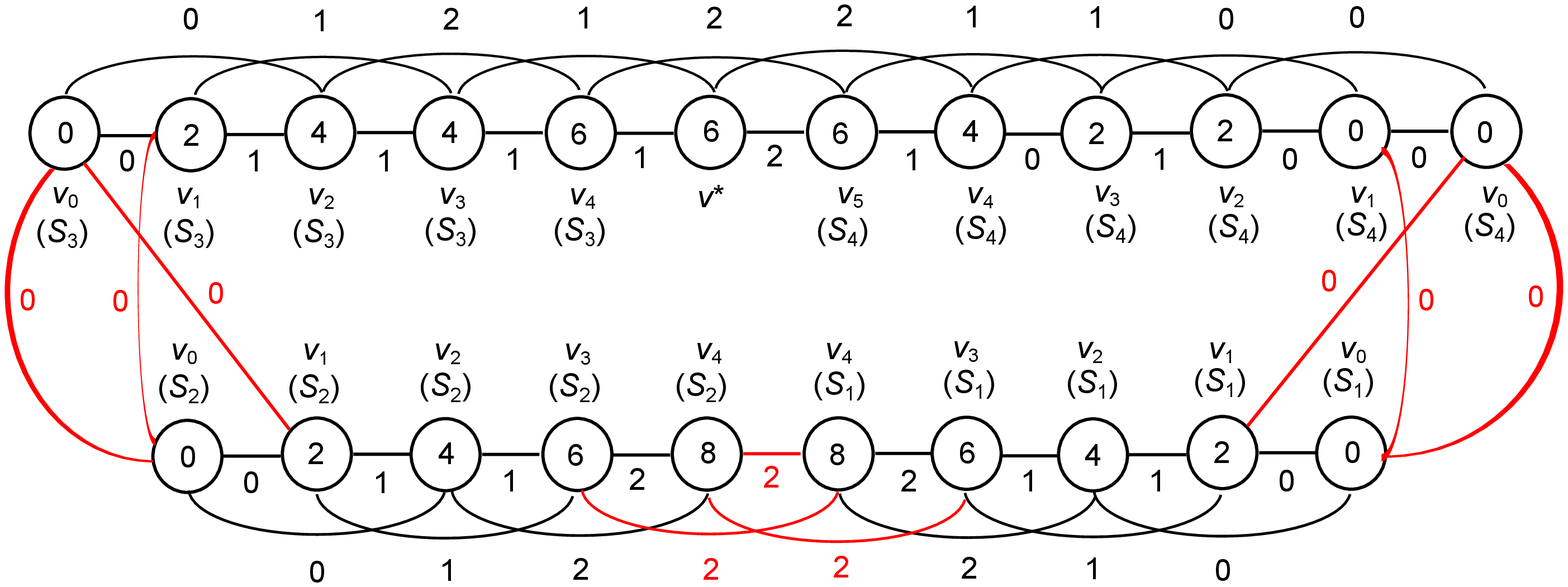}}
    \end{center}
    \caption{Weighting $w^\star$}
    \label{circulant_przyklad_1}
\end{figure}

As $s$ is even, the final labels of edges $e$ for which $w^\star(e)\in\{0,2\}$ we derive using the formula (\ref{nowe_we}), and the remaining ones ($w^\star(e)=1$) - the formula (\ref{nowe_we2}). Two Euler walks - $(v_1(S_1),v_2(S_1),v_3(S_1),v_1(S_1))$ and $(v_1(S_2),v_2(S_2),v_3(S_2),v_1(S_2))$ - have length $3$, so $|V_0|=2$. Let $V_0=\{v_1(S_1),v_1(S_2)\}$, so we increase labels of those vertices by $1$. Final weighting is presented on Figure \ref{circulant_przyklad_1a}.

\begin{figure}[h]
    \begin{center}
    	\scalebox{0.5}{\includegraphics{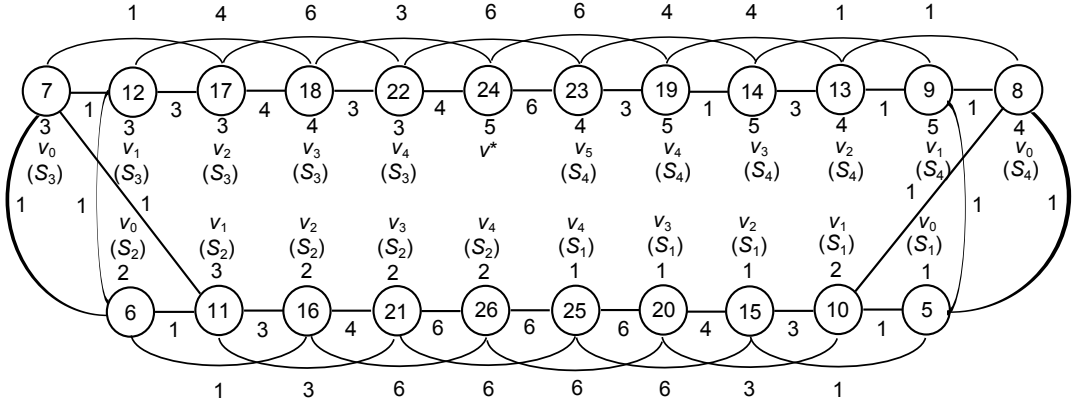}}
    \end{center}
    \caption{Final weighting $w$}
    \label{circulant_przyklad_1a}
\end{figure}
\end{myexample}

\section{Proof of Theorem \ref{theorem_s}}

As in the case of the proof of Theorem \ref{theorem_tvs} we split graph into segments and start with labelling each of the separately. Then we modify the labels in order to obtain desired irregular weighting.

Let us start with two technical lemmas. First of them is in some way extended version of Lemma \ref{lemat_pojedyncze}.

\begin{mylemma}\label{lemat_pojedyncze_bis}{\ }

\noindent{}Let $R=R^{(k)}$ be graph on $2k$ vertices $v_1,v_2,\dots,v_{2k}$ ($k\geq 1$) with edge set consisting of the pairs $(v_i,v_{i+j})$, where $i=1,\dots,2k-1$ and ${j=1,2,\dots,\min\{k,2k-i\}}$. Assume there is a label $l(v_i)$ assigned to every vertex $v_i$, where
\begin{eqnarray*}
l(v_i)=
\begin{cases}
0 &\text{for }  i\leq k\\
2(i-k) &\text{for }  i>k.
\end{cases}
\end{eqnarray*}
Then there exists a weighting $w:E(R)\rightarrow \{0,1,2\}$ such that:
\begin{enumerate}[(i)]
\item
For every vertex $v_i$, $1\leq i \leq 2k$:
\begin{eqnarray*}
\sum_{e\ni v_i}w(e)+l(v_i)=2i.
\end{eqnarray*}
\item
The subgraph of $R$ formed by the edges labelled $1$ and $2$ contains all its vertices and possesses a subgraph $F$ being either a Hamiltonian cycle or a single edge labelled $2$.
\end{enumerate}
\end{mylemma}

\noindent{}\textbf{Proof.}

\noindent{}In order to label the edges we use the variant of Algorithm \ref{algorytm1}.

\begin{myalgorithm}\label{algorytm1_mod}{\ }
\begin{enumerate}[(i)]
\item
Let $i$ be the lowest integer such that all the vertices $v_j$ for $1\leq j\leq i-1$ are closed and let
\begin{eqnarray*}
p=2i-\sum_{j=\max\{1,i-k\}}^{i-1}w((v_j,v_i))-l(v_i).
\end{eqnarray*}
If $1\leq p\leq k$, assign to the edges labels
\begin{eqnarray*}
w((v_i,v_j))=
\begin{cases}
0 &\text{for }  i+1\leq j \leq i+k-p\\
1 &\text{for }  i+k-p+1\leq j \leq i+k.
\end{cases}
\end{eqnarray*}
If $p\geq k+1$, then set
\begin{eqnarray*}
w((v_i,v_j))=
\begin{cases}
1 &\text{for }  i+1\leq j \leq i+2k-p\\
2 &\text{for }  j\geq i+2k-p+1.
\end{cases}
\end{eqnarray*}
Proceed to step $(ii)$.
\item
Close the vertex $v_i$. If $i=2k+1$ then STOP. Otherwise go back to step $(i)$.
\end{enumerate}
\end{myalgorithm}

The proof of part $(i)$ follows directly from the proof of Lemma \ref{lemat_pojedyncze} we only have to omit vertex $v_0$ with all incident edges.

In order to prove part $(ii)$ observe that if $k\geq 3$ is even then the Hamilton cycle $F$ is defined by the vertex sequence
\begin{eqnarray*}
\begin{array}[]{l}
v_{k+1}, v_1, v_k, v_2, v_{k-1},\dots,v_{k/2-1}, v_{k/2+2}, v_{k/2}, v_{k/2+1}, v_{3k/2}, v_{3k/2-1}, v_{3k/2+1}, \\
v_{3k/2-2}, v_{3k/2+2}, \dots, v_{k+3}, v_{2k-3}, v_{k+2}, v_{2k-2}, v_{2k-1}, v_{2k}, v_{k+1},
\end{array}
\end{eqnarray*}
and otherwise by the sequence
\begin{eqnarray*}
\begin{array}[]{l}
v_{k+1}, v_1, v_k, v_2, v_{k-1}, \dots, v_{k/2-3/2}, v_{k/2+5/2}, v_{k/2-1/2}, v_{k/2+3/2}, v_{k/2+1/2}, v_{3k/2-1/2}, \\
v_{3k/2+1/2}, v_{3k/2-3/2}, v_{3k/2+3/2}, \dots, v_{k+3}, v_{2k-3}, v_{k+2}, v_{2k-2}, v_{2k-1}, v_{2k}, v_{k+1}.
\end{array}
\end{eqnarray*}

More precisely, it is formed by the edges
\begin{eqnarray*}
\begin{array}[]{ll}
(v_i, v_{k-i+1}) &\text{for }  1\leq i \leq \lfloor \frac{k}{2} \rfloor,\\
(v_i, v_{k-i+2}) &\text{for }  1\leq i \leq \lceil \frac{k}{2} \rceil,\\
(v_i, v_{3k-i}) &\text{for }  k+2\leq i \leq \lceil \frac{3k}{2} \rceil-1,\\
(v_i, v_{3k-i-1}) &\text{for }  k+2\leq i \leq \lfloor \frac{3k}{2} \rfloor-1,\\
(v_{\lfloor \frac{k}{2} \rfloor+1}, v_{\lfloor \frac{3k}{2} \rfloor}), \\
(v_{k+1}, v_{2k}), \\
(v_{2k-1}, v_{2k}), \\
(v_{2k-2}, v_{2k-1}).
\end{array}
\end{eqnarray*}

If $k=3$, the Hamilton cycle $F$ is formed by the vertex sequence $v_1$, $v_3$, $v_2$, $v_5$, $v_6$, $v_4$, $v_1$, and if $k=2$ - by the sequence $v_1$,$v_3$,$v_4$,$v_2$,$v_1$. If $k=1$, then $(v_1,v_2)$, being the only edge of $R$, is labelled $2$, so $F=R$.\qed

The second lemma guarantees the existence of irregular labellings of graphs $S^{(k)}$ and $R^{(k)}$ with $-1$, $0$ and $1$ such that the resulting weighted degree sequence consists of consecutive integers (with one possible exception).

\begin{mylemma}\label{lemat_jedynki}{\ }

\noindent{}Let $S^{(k)}$ and $R^{(k)}$ be the graphs defined as in the Lemmas \ref{lemat_pojedyncze} and \ref{lemat_pojedyncze_bis}, ${k\geq 2}$. If $G\cong S^{(k)}$ or $G\cong R^{(k)}$, then there exists a labelling ${f:E(G)\rightarrow \{-1,0,1\}}$ such that the weighted degrees of vertices of $G$ are distinct integers from the set $\{0,1,2,\dots,|V(G)|-1\}$ when $k$ is even and distinct integers from the set $\{-1,1,2,\dots,|V(G)|-1\}$ when $k$ is odd.
\end{mylemma}

\noindent{}\textbf{Proof.}

\noindent{}We label the edges using the formula:
\begin{eqnarray}\label{wzor_jedynki}
f(v_i,v_j)=\left\{
\begin{array}[]{ll}
1 &\text{for }  2k+1-|V(G)|\leq i\leq k, 1\leq j\leq k,\\
-1&\text{for }  k+1\leq i \leq 2\left\lfloor\frac{k}{2}\right\rfloor+k-3, j=i+1,\\
-1 &\text{for }  i=2k-1, j=2k.
\end{array}
\right.
\end{eqnarray}
If $k$ is odd, we also set $f(v_{2k-2},v_{2k})=-1$. All remaining edges obtain label $0$.

Given weighting fulfils the conclusion of the Lemma, as
\begin{eqnarray*}
wd(v_i)=
\begin{cases}
|V(G)|-k+i-1 &\text{for }  i\leq k,\\
2k-i &\text{for }  k+1\leq i\leq 2k-1,\\
\end{cases}
\end{eqnarray*}

and
\begin{eqnarray*}
wd(v_{2k})=
\begin{cases}
0 &\text{for }  k\bmod 2=0,\\
-1 &\text{for }  k\bmod 2=1.\\
\end{cases}
\end{eqnarray*}\qed

Let us go back to the proof of Theorem \ref{theorem_s}. As in the proof of Theorem \ref{theorem_tvs}, the construction looks different for distinct relations between $n$ and $k$, so we consider few cases.

\vspace{10pt}
\noindent{}\textit{Case $1$: $n=2k+1$.}

In this case $C_n^k$ is isomorphic to $K_n$ and the proof of the equality
\begin{eqnarray*}
s(C_n^k)=s(K_n)=3
\end{eqnarray*}
may be found e.g. in G. Chartrand et al. \cite{ref_ChaJacLehOelRuiSab1}.

\vspace{10pt}

\noindent{}\textit{Case $2$: $n=2k+2$.}

We use the same edge weighting $w:E(C_n^k)\rightarrow \{1,2\}$ as in the proof of Theorem \ref{theorem_tvs} (Case $2$). Let us remind that the sums of weights of edges incident with vertices $v_i,-\lceil\frac{k}{2}\rceil\leq i \leq \lfloor\frac{k}{2}\rfloor$ are distinct integers from the set $\{2k,\dots,3k\}$, and the sums of edges incident with remaining vertices - distinct integers from the set $\{3k,\dots,4k\}$.

We add $1$ to the labels of edges $(v_i,v_{k+1})$ for $-k\leq i \leq -\lceil\frac{k}{2}\rceil-1$ or ${\lfloor\frac{k}{2}\rfloor+1 \leq i \leq k}$. This way all the weighted degrees become distinct integers from the set $\{2k,\dots,4k\}\cup\{5k\}$. As we use only labels $1$, $2$ and $3$, irregularity strength equals $s(C_{2k+2}^k)=3$.

\vspace{10pt}

\noindent{}\textit{Case $3$: $n>2k+2$.}

We can express the order of $G=C_n^k$ as $n=4kt+r$, where $t\geq 0$ and $1\leq r\leq 4k$ are some integers.

We split $G$ into some number of segments, each except two isomorphic to $R^{(k)}$, one to $R^{(g)}$ and one to $R^{(h)}$, where $g,h\leq k, |g-h|\leq 1$. For some values of $r$ it is necessary to include one additional vertex or a copy of $R^{(r^\prime)}$ or $S^{(r^\prime)}$, where $r^\prime=\left\lfloor\frac{r}{2}\right\rfloor$. Using Lemmas \ref{lemat_pojedyncze_bis} and \ref{lemat_jedynki} we obtain the labellings of segments, then we expand it to all the edges of $G$.

The construction depends on the exact relation between $r$ and $k$, so we consider two cases. For convenience we analyse first the case when $2k+2\leq r\leq 4k$, and then $1\leq r\leq 2k+1$.

\vspace{10pt}
\noindent{}\textit{Case $3.1$: $2k+2\leq r\leq 4k$.}

Let $g=\left\lfloor\frac{r}{4}\right\rfloor$. Consider $2t$ graphs $S_1, S_2,\dots S_{2t}$ isomorphic to $R^{(k)}$ and two graphs $S_{2t+1}\cong R^{(h)}$ and $S_{2t+2}\cong R^{(g)}$, where $h=g$ if ${r \bmod 4 \in \{0,1\}}$ and $h=g+1$ if ${r \bmod 4 \in \{2,3\}}$.

Denote the vertices of $G$ by $v_0(G),\dots,v_{n-1}(G)$. We split $G$ into segments in such a way that the vertices of graphs $S_j, 1\leq j\leq 2t$ (denoted $v_1(S_j),\dots,v_{2k+1}(S_j)$), $S_{2t+1}$ (denoted $v_1(S_{2t+1}),\dots,v_{2h}(S_{2t+1})$) and $S_{2t+2}$ (denoted $v_1(S_{2t+2}),\dots,v_{2g}(S_{2t+2})$) are identified with consecutive subsets of vertices of $G$. For increasing indices of vertices of $G$ the indices in $S_j$ increase when $j=2i-1$ and decrease when $j=2i$, $1\leq i\leq t+1$. In other words, if the vertices $v_i(G)$ and $v_{i+1}(G)$ (or $v_{n-1}(G)$ and $v_{0}(G)$) belong to two neighbouring segments $S_j$ and $S_{j+1}$ (respectively $S_{2t+2}$ and $S_1$), then either they both have maximum indices in the segments (i.e. $v_i(G)=v_{2k}(S_j)$ and $v_{i+1}(G)=v_{2k}(S_{j+1})$, or $v_i(G)=v_{2h}(S_j)$ and $v_{i+1}(G)=v_{2g}(S_{j+1})$), or $v_i(G)=v_{1}(S_j)$ and $v_{i+1}(G)=v_{1}(S_{j+1})$. If $r \bmod 4 \in \{1,3\}$, we insert additional vertex $v^\star=v_{n-2h-1}(G)$ between $S_{2t+1}$ and $S_{2t+2}$.

We label $G$ in two steps. We begin with temporary weighting $w^\star:E(G)\rightarrow \{0,1,2\}$, next we modify the labels in order to obtain final weighting $w:E(G)\rightarrow \{1,\dots,s\}$.

The construction of $w^\star$ begins with labelling the segments using Algorithm \ref{algorytm1_mod}. Then we label remaining edges of $G$ (i.e. $e\in E(G)\setminus \bigcup_{j=1}^{2t+2}E(S_j)$).

If two consecutive segments are $S_{2i-1}$ and $S_{2i}$, $1\leq i\leq t$ (so two consecutive vertices belonging to them have degree $2k$), we put $2$ on all the edges joining these segments. The edges joining $S_{2i}$ with $S_{2i+1}$, $1\leq i\leq t$, and $S_{2t+2}$ with $S_1$ obtain label $0$.

Observe that the weighted degrees do not change as the sums of weights of edges joining different segments of $G$ are equal to $l(v)$ from Lemma \ref{lemat_pojedyncze_bis}.

If $r\bmod 2=0$, we put $2$ on the edges
\begin{eqnarray*}
\begin{array}[]{ll}
(v_{h+i}(S_{2t+1}),v_{2g-j+1}(S_{2t+2})) &\text{for } 1\leq i\leq h, 1\leq j\leq i,
\end{array}
\end{eqnarray*}
and otherwise on the edges
\begin{eqnarray*}
\begin{array}[]{ll}
(v_{h+i}(S_{2t+1}),v_{2g-j+2}(S_{2t+2})) &\text{for } 2\leq i\leq h, 1\leq j\leq i-1,\\
(v^\star,v_{h+i}(S_{2t+1})) &\text{for } 1\leq i\leq g+1,\\
(v^\star,v_{2g-i+1}(S_{2t+2})) &\text{for } 1\leq i\leq g+1.
\end{array}
\end{eqnarray*}

Observe that the last operation may increase by $2$ weighted degrees of some vertices $v\in V(S_{2t+2})$, however the conditions
$$wd^\star(v)\leq 2k {\ \ \ }\text{for } v\in V(S_{2t+1})\cup V(S_{2t+2})\cup \{v^\star\}$$
and
$$|\{v\in V(S_{2t+1})\cup V(S_{2t+2})\cup \{v^\star\}: wd^\star(v)=2i\}|\leq 2{\ \ \ }\text{for }i=1,\dots,k$$
still hold as $wd(v^\star)\neq wd(v_{2g}(S_{2t+2}))$.

In order to finish the construction of $w^\star$, we assign $0$ to all not labelled edges of $G$.

Finally we obtain weighting $w^\star:E(G)\rightarrow \{0,1,2\}$ such that all weighted degrees are even, segments $S_j$ contain subgraphs $F$ defined in Lemma \ref{lemat_pojedyncze_bis}, and moreover
$$|\{v\in V(G): wd^\star(v)=2i\}|\leq 2t+2{\ \ \ }\text{for }i=1,\dots,k.$$

Let $s=\lceil \frac{n+2k-1}{2k} \rceil$ (observe that $s$ is odd). We define new weighting $w:E(G)\rightarrow \{0,\frac{s+1}{2}, s\}$ as follows:
\begin{eqnarray}\label{nowe_we_mod}
w(e)=\frac{s-1}{2}w^\star(e)+1{\ \ \ }\text{for } e\in E(G).
\end{eqnarray}

After this modification weighted degrees in any segment $S_j$ differ by at least $s-1$. As we can not label vertices, we are going to distinguish their degrees by changing labels on the edges of subgraphs $F$ defined in Lemma \ref{lemat_pojedyncze_bis}.

We decrease the degrees of vertices $v\in V(S_j), 1\leq j\leq 2t$, in the following way. If $j\bmod 2=0$, we decrease the weight of every edge in $F$ by $\frac{j}{2}$. Otherwise we decrease it alternately by $\frac{j-1}{2}$ and $\frac{j+1}{2}$ (each subgraph of $F$ is a cycle of  even length).

If $h>1$, we decrease the labels in $S_{2t+1}$ alternately by $t$ and ${t+1}$. Otherwise we decrease the weight of the only edge in this segment $(v_1(S_{2t+1}),v_2(S_{2t+1}))$ by $2t+1$.

If $g=k$ (it is possible only if $r=4k$), we decrease in the same way the labels in $S_{2t+2}$. Otherwise we do not change neither the weights in $S_{2t+2}$ nor the weights of edges incident with $v^\star$ (if this vertex is considered).

This way we decrease weighted degrees in $S_{j}$ by distinct integers from the set $\{0,1,\dots,2t+1\}$ or $\{1,2,\dots,2t+2\}$ and finally distinguish all the weighted degrees in $G$, as $2t+2\leq s-1$.

\vspace{10pt}
\noindent{}\textit{Case $3.2$: $1\leq r\leq 2k+1$.}

We start with finding an irregular weighting of $C_{4kt}^k$ using the method described in previous paragraph. Let $s=2t+2$. Observe that the maximum weight used so far is $s-1$ and two highest weighted degrees are equal $wd(v^\star)=2k(s-1)-1$ and $wd(v^{\star\star})=2k(s-1)-2$ for vertices $v^\star=v_{2k}(S_1)$ and $v^{\star\star}=v_{2k}(S_2)$. Observe also that all the edges between $S_1$ and $S_2$ are labelled $s-1$.

If $r=1$, we insert additional vertex $v_0$ between $v^{\star}$ and $v^{\star\star}$ removing and adding the edges in order to obtain $C_{4kt+1}^k$. We put $s-1$ on all the edges incident with $v_0$. The weighted degree of $v_0$ reaches the value $wd(v_0)=2k(s-1)$ and the remaining weighted degrees do not change, so this way we obtain irregular weighting of $G$, as $\left\lceil \frac{4kt+1+2k-1}{2k}\right\rceil=s-1$. 

If $r=2$, we insert two vertices $v_1$ and $v_2$ between $v^{\star}$ and $v^{\star\star}$, and put $s-1$ on all new edges except $(v_1,v^{\star})$ and $(v_1,v_2)$, that are labelled $s$. This way the weighted degrees of vertices $v\in V(C_n^k)\setminus\{v^{\star},v_1,v_2\}$ do not change, so they are distinct integers from the set 
$$\{2k,\dots,2k(s-1)-2\}.$$
Moreover
\begin{eqnarray*}
\begin{array}[]{l}
wd(v^{\star})=2k(s-1),\\
wd(v_1)=2k(s-1)+2,\\
wd(v_2)=2k(s-1)+1.
\end{array}
\end{eqnarray*}
From the above it follows that $w$ is irregular weighting of $G$. 

If $r=3$, we insert three vertices $v_1$, $v_2$ and $v_3$ between $v^{\star}$ and $v^{\star\star}$, and next put $s-1$ on all new edges except the following:
$$w(v^{\star},v_1)=s-2,$$
$$w(v^{\star\star},v_2)=w(v_1,v_2)=w(v_1,v_3)=w(v_2,v_3)=s.$$
The resulting weighted degrees of vertices from the set
$$v\in V(C_n^k)\setminus\{v^{\star},v^{\star\star},v_1,v_2, v_3\}$$
are distinct integers from the set
$$\{2k,\dots,2k(s-1)-3\}.$$
Moreover
\begin{eqnarray*}
\begin{array}[]{l}
wd(v^{\star})=2k(s-1)-2,\\
wd(v^{\star\star})=2k(s-1)-1,\\
wd(v_1)=2k(s-1)+1,\\
wd(v_2)=2k(s-1)+3,\\
wd(v_3)=2k(s-1)+2,
\end{array}
\end{eqnarray*}
so $w$ is irregular weighting of $G$. 

Let us move to the case when $4\leq r\leq 2k+1$. Observe that this time $\left\lceil \frac{4kt+r+2k-1}{2k} \right\rceil=s$.

We insert between $v^{\star}$ and $v^{\star\star}$ the graph $H$, where
\begin{eqnarray*}
H\cong
\begin{cases}
S^{(r^\prime)} &\text{for } r\bmod 2=1,\\
R^{(r^\prime)} &\text{for } r\bmod 2=0, 
\end{cases}
\end{eqnarray*}
and
$$r^\prime=\left\lfloor\frac{r}{2}\right\rfloor.$$

Let $V(H)=\{v_0(H),v_1(H),\dots,v_{|V(H)|-1}(H)\}$. We insert $H$ between $v^{\star}$ and $v^{\star\star}$ in such a way that $v^{\star}$ is adjacent to
$$v_j(H) {\ \ \ }\text{for }0\leq j\leq\min\{k-1,|V(H)|-1\},$$
and $v^{\star\star}$ to
$$v_{|V(H)|-j-1}(H) {\ \ \ }\text{for }0\leq j\leq\min\{k-1,|V(H)|-1\}.$$

We put $s-1$ on all new edges, including the edges of $H$. The resulting weighting $w_1:E(G)\rightarrow \{1,2,\dots,s-1\}$ assigns the weighted degree $2k(s-1)$ to every vertex of $H$ and does not change the degrees of the remaining vertices of $G$ (including $v^{\star}$ and $v^{\star\star}$). In order to distinguish the weighted degrees of vertices $v\in V(H)$ we put the labels $f$ on the edges of $H$ using the formula (\ref{wzor_jedynki}) (see Lemma \ref{lemat_jedynki}). If $k\bmod 2=0$, we define $w:E(G))\rightarrow \{1,2,\dots,s\}$ using the formula:

\begin{eqnarray}\label{wzor_nowe_s_z_jedynkami}
w(e)=\left\{
\begin{array}[]{ll}
w_1(e) &\text{for } e\in E(G))\setminus E(H),\\
w_1(e)+f(e) &\text{for } e\in E(H).
\end{array}
\right.
\end{eqnarray}

As $s\geq 3$, new edge labels fulfil the condition $1\leq w(e) \leq s$ and all the weighted degrees are distinct, as according to the Lemma \ref{lemat_jedynki} the weighted degrees of vertices $v\in V(H)$ become distinct integers from the set
$$\{2k(s-1),2k(s-1)+1,\dots,2k(s-1)+r\},$$
and the degrees of vertices $v\in V(G)\setminus V(H)$ obtain distinct values from the set
$$\{2k,2k+1,\dots,2k(s-1)-1\}$$
(they do not change after inserting the graph $H$).

If $k\bmod 2=1$, we also use the formula (\ref{wzor_nowe_s_z_jedynkami}) in order to define new weighting $w$. However this time
$$wd(v^{\star})=wd(v_{2r\prime})=2k(s-1)-1.$$
If $r\leq 2k$, we change the edge label
$$w(v^{\star},v_{r\prime})=s.$$

As it can be easily observed, after this operation
\begin{eqnarray}
\begin{array}[]{l}
wd(v_{r\prime})=2k(s-1)+r+1,\\
wd(v^{\star})=2k(s-1),\\
wd(v_{2r\prime})=2k(s-1)-1,
\end{array}
\end{eqnarray}
so $w$ is the desired irregular weighting of $G$, as the weighted degrees of vertices $v\in V(C_n^k)\setminus(V(G)\cup \{v^{\star}\})$ are distinct and lower than $2k(s-1)-1$, and the weighted degrees of vertices $v\in V(G)\setminus\{v_{r\prime},v_{2r\prime}\}$ are distinct integers from the set $\{2k(s-1)+2,\dots,2k(s-1)+r-1\}$.

If $r= 2k+1$, then  $(v^{\star},v_{r\prime})\not\in E(G)$, so we have to proceed in an different way. We redefine the weighting assigning
$$w(v_{k}(H),v_{2k}(H))=s+1.$$

This way $wd(v_{k}(H))$ and $wd(v_{2k}(H))$ increase by $1$, so the weighted degrees of vertices $v\in V(G)$ reach distinct values. Observe that $w(v_{k}(H),v_{2k}(H))=s+1=\lceil \frac{n+2k-1}{2k} \rceil+1$.

Assume that the labels $1,2,\dots,s$ are enough to construct an irregular weighting when $r=2k+1$ and $k\bmod 2=1$. In such a case the weighted degrees of vertices $v\in V(G)$ have to reach all the values of the set $\{2k,\dots,4kt+4k\}$, what means that their sum is equal to the odd number $k(2t+3)(4kt+2k+1)$. However, it is impossible as the sum of degrees has to be even. The contradiction proves that the use of label $s+1$ is necessary.\qed

\vspace{10pt}

\begin{myexample}[Irregular labelling of $C_{19}^3$]{\ }
\\ \noindent{}We have $n=19$, $k=3$ $t=1$, $r=7=2k+1$, so first we find the labelling of $C_{12}^3$. Initial labelling with $0$, $1$ and $2$ are presented on the Figure \ref{circulant_przyklad_2}. The edges joining different segments are red.

\begin{figure}[h]
    \begin{center}
    	\scalebox{0.5}{\includegraphics{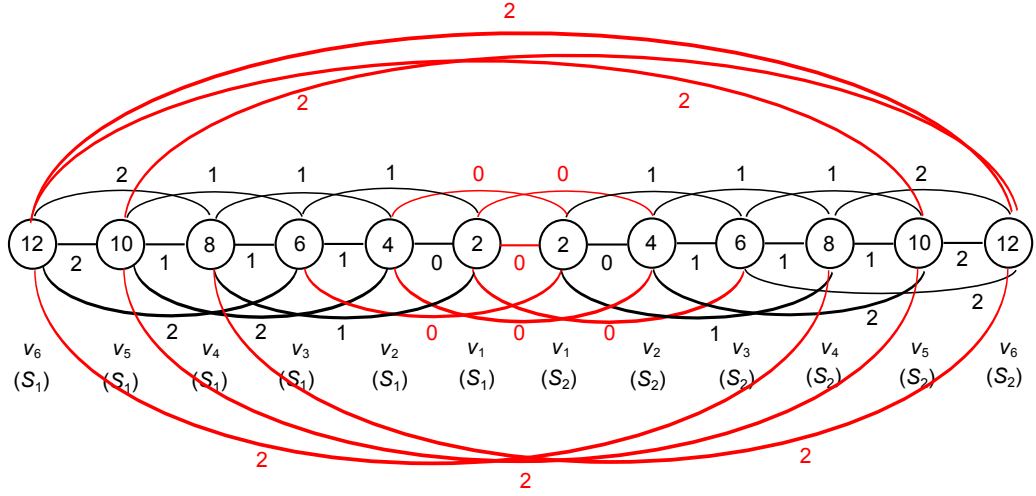}}
    \end{center}
    \caption{Initial labelling of the edges of $C_{12}^3$ with $0$, $1$ and $2$}
    \label{circulant_przyklad_2}
\end{figure}

Observe that $s(C_{12}^3)=s-1=3$ is odd, so we modify the edges using formula (\ref{nowe_we}). Next we modify the weights of edges of subgraphs $F$, decreasing by $1$ half of the labels in $S_1$ all the labels in $S_2$. Final weighting of $C_{12}^3$ is presented in the Figure \ref{circulant_przyklad_2a} (both cycles are red).

\begin{figure}[h]
    \begin{center}
    	\scalebox{0.5}{\includegraphics{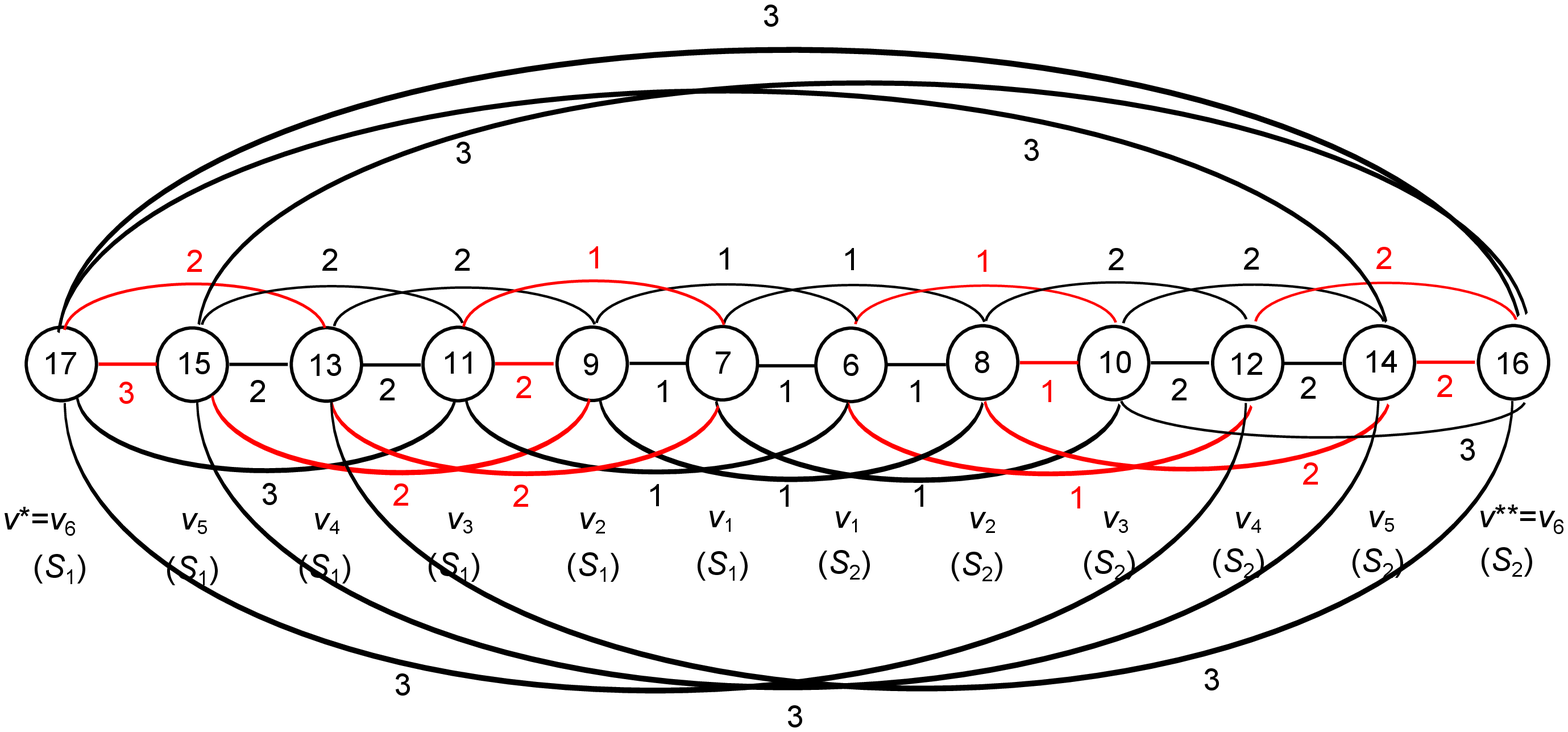}}
    \end{center}
    \caption{Final weighting of $C_{12}^3$}
    \label{circulant_przyklad_2a}
\end{figure}

\newpage
In the next step we construct weighting $f$ of graph $H=S^{(3)}$ as in the Lemma \ref{lemat_jedynki}. It is presented on the Figure \ref{circulant_przyklad_2b}.

\begin{figure}[h]
    \begin{center}
    	\scalebox{0.5}{\includegraphics{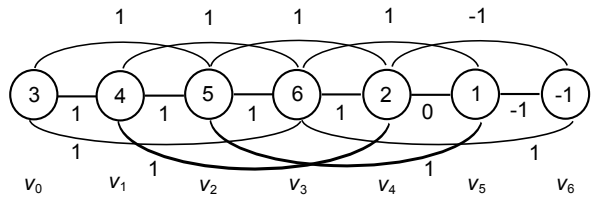}}
    \end{center}
    \caption{Lemma \ref{lemat_jedynki} - weighting $f$ of $H=S^{(3)}$}
    \label{circulant_przyklad_2b}
\end{figure}

In order to finish, we insert graph $H\cong S^{(3)}$ between the vertices $v^\star=v_6(S_1)$ and $v^{\star\star}=v_6(S_2)$, putting $s-1$ on all new edges $e\not\in E(H)$, and $s-1+f(e)$ on all the edges $e\in E(H)$. As $r=2k+1$, we assign $w(v_3(H),v_6(H))=s-1+f(v_3(H),v_6(H))+1=5$. This way we obtain final irregular weighting as on the Figure \ref{circulant_przyklad_2c}.

\begin{figure}[h]
    \begin{center}
    	\scalebox{0.5}{\includegraphics{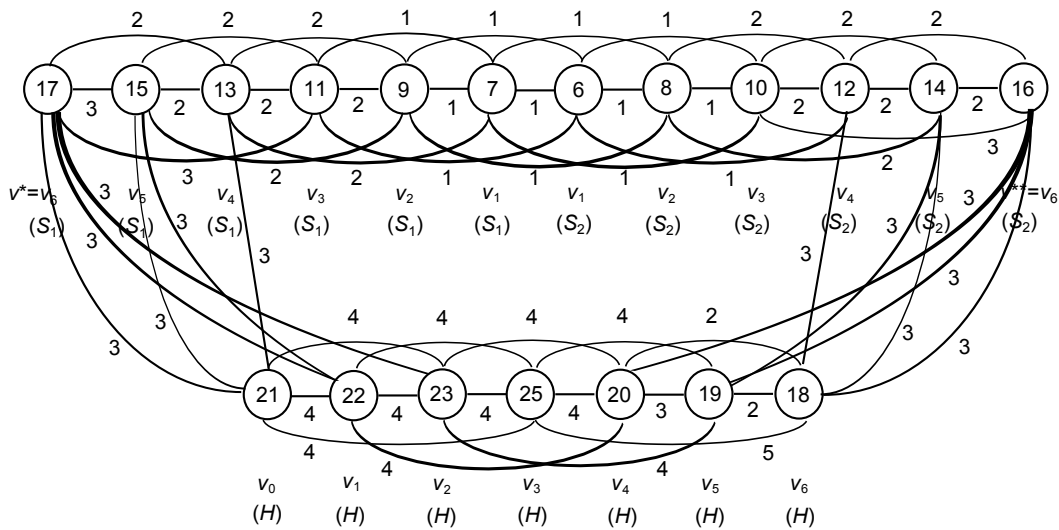}}
    \end{center}
    \caption{Final weighting of $C_{19}^3$}
    \label{circulant_przyklad_2c}
\end{figure}

\end{myexample}

\end {document}